\documentclass[11pt]{article}
\usepackage{latexsym,amssymb}
\def\demo{\noindent{\bf Proof. }}
\def\QED{\hfill$\Box$}
\newtheorem{Theorem}{Theorem}[section]
\newtheorem{Lemma}[Theorem]{Lemma}
\newtheorem{Corollary}[Theorem]{Corollary}

\newtheorem{Remark}[Theorem]{Remark}

\newtheorem{Definition}[Theorem]{Definition}

\begin{document}
\topmargin3mm
\hoffset=-1cm
\voffset=-1.5cm
\

\medskip

\begin{center}
{\large\bf Edge ideals of clique clutters of comparability graphs 
and the normality of monomial ideals
}
\vspace{6mm}\\
\footnotetext{2000 {\it Mathematics Subject 
Classification\/}. Primary 13A30; Secondary 52B20, 06A07,13B22.} 
\footnotetext{{\it Key words and phrases\/}.  
normal ideal, normally torsion free, max-flow min-cut, edge ideal, 
comparability graph, integer rounding property,
clique clutter, poset, integer decomposition property.}

\medskip

Luis A. Dupont and Rafael H. Villarreal\footnote{Partially supported by CONACyT 
grant 49251-F and SNI.} 
\\ 
{\small Departamento de Matem\'aticas}\vspace{-1mm}\\ 
{\small Centro de Investigaci\'on y de Estudios Avanzados del
IPN}\vspace{-1mm}\\   
{\small Apartado Postal 14--740}\vspace{-1mm}\\ 
{\small 07000 M\'exico City, D.F.}\vspace{-1mm}\\ 
{\small e-mail: {\tt vila@math.cinvestav.mx}}\vspace{4mm}
\end{center}
\date{}

\begin{abstract} 
\noindent The normality of a monomial ideal is 
expressed in terms of lattice points of blocking polyhedra and the integer 
decomposition property. For edge ideals of clutters this 
property characterizes normality. Let $G$ be the comparability graph of
a finite poset. If ${\rm cl}(G)$ is the clutter of maximal
cliques of $G$, we prove that ${\rm cl}(G)$ satisfies the max-flow 
min-cut property and that its edge ideal is normally torsion free.
Then we prove that edge ideals of complete admissible uniform clutters are 
normally torsion free. 
\end{abstract}

\section{Introduction}

Let $R=K[x_1,\ldots,x_n]$ be a polynomial ring 
over a field $K$ and let $I$ be a monomial ideal of $R$. We are
interested in determining what families of monomial ideals have the 
property that $I$ is {\it normal\/} or {\it normally torsion free}. An
aim here is to explain how these two algebraic properties interact with
combinatorial optimization and linear programming problems. 
Recall that $I$ is called normal
(resp. normally torsion free ) if $I^i=\overline{I^i}$ (resp.
$I^i=I^{(i)}$) for all $i\geq 1$, where $\overline{I^i}$ and 
$I^{(i)}$ denote the integral closure of the 
$i${\it th} power of $I$ and the $i${\it th}
symbolic power of $I$ 
respectively  (see the beginning of 
Sections~\ref{normali-rounding} and \ref{ntf-normal-sec} for the
precise definitions of $\overline{I^i}$ and $I^{(i)}$). 
If $\overline{I}=I$, the ideal $I$ is called {\it integrally
closed\/}.

The contents of this paper are as follows. In
Section~\ref{normali-rounding} 
we study the normality
of monomial ideals. We are able to characterize this property
in terms of blocking polyhedra and the integer decomposition property
(see Theorem~\ref{normal-iff-integerdec}). For integrally closed
ideals this property characterizes normality 
(see Corollary~\ref{normal-iff-integerdec-ic}). 
As a consequence using a
result of Baum and Trotter \cite{baum-trotter} we describe the
normality of a monomial ideal in terms of the integer rounding property 
(see Corollary~\ref{crit-rounding-normali}). 

Before introducing the main results of Sections~\ref{maximal-cliques}
and \ref{ntf-normal-sec}, let us recall some notions that will play an
important role in what follows. Let $\mathcal{C}$ be a {\it
clutter\/} with finite vertex set  
$X=\{x_1,\ldots,x_n\}$, 
that is, $\mathcal{C}$ is a family of subsets of $X$, called edges, none 
of which is included in another. The set of
vertices and edges of $\mathcal{C}$ are denoted by $V(\mathcal{C})$
and $E(\mathcal{C})$ respectively. The {\it incidence matrix\/} of
$\mathcal{C}$ is the vertex-edge matrix whose columns are the 
characteristic vectors of the edges of $\mathcal{C}$. 
The {\it edge ideal\/} of $\mathcal{C}$, 
denoted by $I(\mathcal{C})$, is the ideal of $R$
generated by all monomials $\prod_{x_i\in e}x_i$ such 
that $e\in E(\mathcal{C})$. 

Let $P=(X,\prec)$ be a {\it partially ordered set\/} ({\it poset} for
short) on the finite vertex set $X$ 
and let $G$ be its {\it comparability graph\/}. Recall that the vertex 
set of $G$ is $X$ and the edge set of $G$ is the set of all unordered
pairs $\{x_i,x_j\}$ such that $x_i$ and $x_j$ are comparable. 
A {\it clique\/} of $G$ is a subset of the set
of vertices of $G$ that induces 
a complete subgraph. The {\it clique clutter\/} of $G$, 
denoted by ${\rm cl}(G)$, is the clutter with vertex set $X$ whose
edges are exactly the maximal cliques of $G$ (maximal with respect to
inclusion).

Our main algebraic result is presented in
Section~\ref{ntf-normal-sec}. It shows that the edge ideal 
$I=I({\rm cl}(G))$ of ${\rm cl}(G)$ is normally torsion free (see 
Theorem~\ref{compara-ntf-normal}). To prove
this result we first show that  
the clique clutter of $G$ has the max-flow min-cut property (see 
Theorem~\ref{comp-max-cliques-mfmc}). Then we 
use a remarkable result of  \cite{clutters} showing that an edge ideal
$I(\mathcal{C})$, of a clutter $\mathcal{C}$, is normally torsion free
if and only if $\mathcal{C}$ 
has the max-flow min-cut property. As an application, we prove
that edge ideals of  
complete admissible uniform clutters are normally torsion free (see
Theorem~\ref{cauc-ntf}). This interesting family of clutters was 
introduced and studied in \cite{floystad}.

Along the paper we introduce most of the 
notions that are relevant for our purposes. Our main references for
combinatorial optimization and commutative 
algebra are \cite{BHer,Schr2,bookthree,monalg}. In these references
the reader will find the undefined terminology and notation
that we use in what follows. 

\section{Normality of monomial ideals}\label{normali-rounding}

Let $R=K[x_1,\ldots,x_n]$ be a polynomial ring over a field $K$, 
let $I$ be a monomial ideal of $R$ generated by
$x^{v_1},\ldots,x^{v_q}$, and let $A$ be the $n\times q$ 
matrix with column vectors $v_1,\ldots,v_q$. As usual, we will use  $x^a$ as
an abbreviation for  
$x_1^{a_1} \cdots x_n^{a_n}$, where $a=(a_i)$ is a vector in 
$\mathbb{N}^n$. Recall that the {\it integral 
closure\/} of $I^i$, denoted by
$\overline{I^i}$, is the ideal of $R$ given by  
\begin{equation}\label{oct12-08}
\overline{I^i}=(\{x^a\in R\vert\, \exists\,
p\in\mathbb{N}\setminus\{0\};(x^a)^{p}\in
I^{pi}\}),
\end{equation}
see for instance \cite[Proposition~7.3.3]{monalg}. The ideal $I$ 
is called {\it normal\/} if $I^i=\overline{I^i}$ for $i\geq 1$. 
In this section we give a characterization of the normality of $I$ in
terms of lattice points of blocking polyhedra. The polyhedron 
$$
B(Q)=\{z\in\mathbb{R}^n\vert\, z\geq 0;\, \langle z,x\rangle\geq 1\mbox{ for all
}x\mbox{ in }Q \}
$$
is called the {\it blocking polyhedron\/} of $Q=Q(A)=\{x\vert\,x\geq 0;\,
xA\geq \mathbf{1}\}$. The polyhedron $B(Q)$ 
is said to have the {\it integer decomposition property\/} if for each
natural number $k$ and for each integer vector $a$ in $kB(Q)$, $a$ is the sum of
$k$ integer vectors in $B(Q)$; see \cite[pp.~66--82]{Schr2}. 

\begin{Theorem}\label{normal-iff-integerdec} The ideal $I$ is normal if and only if 
the blocking polyhedron $B(Q)$ of $Q=Q(A)$ has the integer
decomposition property and all minimal
integer vectors of $B(Q)$ are columns of $A$ 
$($minimal with respect to $\leq$$)$.
\end{Theorem}

\demo First we show the equality $B(Q)=\mathbb{R}_+^n+{\rm
conv}(v_1,\ldots,v_q)$. The right hand side is clearly contained in
the left hand side. Conversely take $z$ in $B(Q)$, then 
$\langle z,x\rangle\geq 1$ for all $x\in Q(A)$ and $z\geq 0$. Let
$\ell_1,\ldots,\ell_r$ be the vertex set of $Q(A)$. In particular 
$\langle z,\ell_i\rangle\geq 1$ for all $i$. Then 
$\langle(z,1),(\ell_i,-1)\rangle\geq 0$ for all $i$. From
\cite[Theorem~3.2]{clutters} 
we get that $(z,1)$ belongs to the cone generated by  
$$
\mathcal{A}'=\{e_1,\ldots,e_n,(v_1,1),\ldots,(v_q,1)\}.
$$ 
Thus $z$ is
in $\mathbb{R}_+^n+{\rm conv}(v_1,\ldots,v_q)$. This completes the
proof of the asserted equality. Hence $B(Q)\cap \mathbb{Q}^n=\mathbb{Q}_+^n+{\rm
conv}_\mathbb{Q}(v_1,\ldots,v_q)$ because the polyhedron $B(Q)$ is
rational. Using this equality and the description of the integral
closure given in Eq.~(\ref{oct12-08}), we readily obtain the equality
\begin{equation}\label{ic-kth-power}
\overline{I^k}=(\{x^a\vert\,
a\in{kB(Q)}\cap\mathbb{Z}^n\})
\end{equation}
for $0\neq k\in\mathbb{N}$. Assume that $I$ is normal, i.e.,
$\overline{I^k}=I^k$ for $k\geq 1$. Let $a$ be an integer vector in
$kB(Q)$. Then $x^a\in I^k$ and consequently $a$ is the sum of $k$
integer vectors in $B(Q)$, that is, $B(Q)$ has the integer
decomposition property. Take a minimal integer vector $a$ in $B(Q)$. 
Then $x^a\in\overline{I}=I$ and we can write $a=\delta+v_i$ for some
$v_i$ and for some $\delta\in\mathbb{N}^n$. Thus $a=v_i$ by the
minimality of $a$. Conversely assume that $B(Q)$ has the integer
decomposition property and all minimal integer vectors of $B(Q)$ are
columns of $A$. Take $x^a\in\overline{I^k}$, i.e., $a$ is an integer
vector of $kB(Q)$. Hence $a$ is the sum of 
$k$ integer vectors $\alpha_1,\ldots,\alpha_k$ in $B(Q)$. Since any
minimal vector of $B(Q)$ is a column of $A$ we may assume that 
$\alpha_i=c_i+v_i$ for $i=1,\ldots,k$. Hence $x^a\in I^k$, as
required. \QED

\begin{Corollary}\label{normal-iff-integerdec-ic} If 
$I=\overline{I}$, then $I$ is normal if and only if 
the blocking polyhedron $B(Q)$ has the integer
decomposition property.
\end{Corollary}

\demo $\Rightarrow$) If $I$ is normal, by Theorem~\ref{normal-iff-integerdec} the
blocking polyhedron $B(Q)$ has the integer decomposition property.

$\Leftarrow$) Take $x^a\in\overline{I^k}$. From
Eq.~(\ref{ic-kth-power}) we get that $a$ is an integer
vector of $kB(Q)$. Hence $a$ is the sum of 
$k$ integer vectors $\alpha_1,\ldots,\alpha_k$ in $B(Q)$. Using 
Eq.~(\ref{ic-kth-power}) with $k=1$, we get that
$\alpha_1,\ldots,\alpha_k$ are in $\overline{I}=I$. Hence $x^a\in I^k$, as
required. \QED

\begin{Corollary}\label{normal-iff-integerdec-edge} If
$I=I(\mathcal{C})$ is the edge ideal of a clutter $\mathcal{C}$, 
then $I$ is normal if and only if the blocking polyhedron $B(Q)$ has
the integer decomposition property. 
\end{Corollary}

\demo Recall that $I$ is an intersection of prime 
ideals (see \cite[Corollary~5.1.5]{monalg}). 
Thus it is seen that 
$\overline{I}=I$. Then the result follows from 
Corollary~\ref{normal-iff-integerdec-ic}. \QED

\begin{Definition}\rm The system $x\geq 0; xA\geq\mathbf{1}$ of linear
inequalities is said to have the {\it integer rounding property\/} 
if 
$${\rm max}\{\langle y,{\mathbf 1}\rangle \vert\, 
y\geq 0; Ay\leq w; y\in\mathbb{N}^q\} 
=\lfloor{\rm max}\{\langle y,{\mathbf 1}\rangle \vert\, y\geq 0;
Ay\leq w\}\rfloor
$$
for each integer vector $w$ for which the right 
hand side is finite. 
\end{Definition}

Systems with the integer rounding property have been widely studied; see
for instance \cite[Chapter~22, pp. 336--338]{Schr}, 
\cite[pp. 82--83]{Schr2}, and the
references there. 

\begin{Corollary}\label{crit-rounding-normali} 
The ideal $I$ is normal ideal if and only if 
the system $xA\geq\mathbf{1};x\geq 0$ has the integer rounding
property. 
\end{Corollary}

\demo According to \cite{baum-trotter} the system $xA\geq\mathbf{1};
x\geq 0$ has the integer rounding property if and only if 
the blocking polyhedron $B(Q)$ of $Q=Q(A)$ has the integer
decomposition property and all minimal
integer vectors of $B(Q)$ are columns of $A$ 
$($minimal with respect to $\leq$$)$ (cf. \cite[p.~82, Eq.~(5.80)]{Schr2}). 
Thus the result follows at once
from Theorem~\ref{normal-iff-integerdec}. \QED

\medskip

There are some other useful characterizations of the normality of a monomial ideal 
\cite[Theorem~4.4]{normali}. 

\section{Maximal cliques of comparability graphs}\label{maximal-cliques}

In
this section we introduce the max-flow min-cut property and prove our
main combinatorial result, 
that is, we prove that the clique clutter of a comparability graph
satisfies the max-flow min-cut property.

\begin{Definition}\rm Let $\mathcal{C}$ be a clutter and let $A$ be
its incidence matrix. The clutter $\cal C$ satisfies the {\it
max-flow min-cut\/}  
property if both sides 
of the LP-duality equation
\begin{equation}\label{jun6-2-03-1}
{\rm min}\{\langle w,x\rangle \vert\, x\geq 0; xA\geq{\mathbf 1}\}=
{\rm max}\{\langle y,{\mathbf 1}\rangle \vert\, y\geq 0; Ay\leq w\} 
\end{equation}
have integer optimum solutions $x$ and $y$ for each non-negative 
integer vector $w$. 
\end{Definition}

Let $\mathcal{C}$ be a clutter. A set of edges of $\cal C$ is {\it
independent\/} or {\it stable} if no two of them have a
common vertex. 
We denote the smallest number of vertices in any 
minimal vertex cover of $\cal C$ by $\alpha_0({\cal C})$ and the 
maximum number of independent edges of ${\cal C}$ by 
$\beta_1({\cal C})$. These two numbers satisfy $\beta_1(\mathcal{C})\leq
\alpha_0(\mathcal{C})$.

\begin{Definition}\label{konig-def}\rm If $\beta_1({\cal
C})=\alpha_0({\cal C})$, we say  that $\cal C$ has the {\it K\"onig
property\/}.
\end{Definition}

Let $\mathcal{C}$ be a clutter on the vertex set
$X=\{x_1,\ldots,x_n\}$ and let $x_i\in X$. Then {\it duplicating\/}
$x_i$ means extending $X$ by a new vertex $x_i'$ and replacing
$E(\mathcal{C})$ by
$$    
E(\mathcal{C})\cup\{(e\setminus\{x_i\})\cup\{x_i'\}\vert\, x_i\in e\in
E(\mathcal{C})\}.
$$
The {\it deletion\/} of $x_i$, denoted by
$\mathcal{C}\setminus\{x_i\}$, is the clutter formed from
$\mathcal{C}$ by deleting the vertex $x_i$ and all edges containing
$x_i$. A clutter obtained from $\mathcal{C}$ by a sequence of deletions and
duplications of vertices is called a {\it parallelization\/}. If 
$w=(w_i)$ is a vector in $\mathbb{N}^n$, we denote by $\mathcal{C}^w$
the clutter obtained from
$\mathcal{C}$ by deleting any vertex $x_i$ with $w_i=0$ and
duplicating $w_i-1$ times any vertex $x_i$ if $w_i\geq 1$. 

The notion of parallelization can be used to give the following characterization
of the max-flow min-cut property which is suitable to study the
clique clutter of the comparability graph of a poset.

\begin{Theorem}{\rm\cite[Chapter~79,
Eq.~(79.1)]{Schr2}}\label{mfmc-iff-cw}  
Let $\mathcal{C}$ be
a clutter. Then $\mathcal{C}$ 
satisfies the max-flow min-cut property if and only if 
$\beta_1(\mathcal{C}^w)=\alpha_0(\mathcal{C}^w)$ for all
$w\in\mathbb{N}^n$.
\end{Theorem}

\begin{Lemma}\label{dup-comm-cliques} Let ${\rm cl}(G)$ be the clutter of 
maximal cliques of a graph $G$. If $G^1$  $($resp. ${\rm cl}(G)^1$$)$ is the graph $($resp.
clutter$)$ obtained from $G$ $($resp. ${\rm cl}(G)$$)$ by duplicating the
vertex $x_1$, then ${\rm cl}(G)^1={\rm cl}(G^1)$.
\end{Lemma}
\demo Let $y_1$ be the duplication of $x_1$. Set $\mathcal{C}={\rm
cl}(G)$. First we prove that $E(\mathcal{C}^1)\subset E({\rm
cl}(G^1))$. Take $e\in E(\mathcal{C}^1$). Case (i): 
Assume $y_1\notin e$. Then $e\in E(\mathcal{C})$. Clearly $e$ is a
clique of $G^1$. If $e\notin E({\rm
cl}(G^1))$, then $e$ can be extended to a maximal clique of $G^1$.
Hence $e\cup\{y_1\}$ must be a clique of $G^1$. Note that
$x_1\notin e$ because $\{x_1,y_1\}$ is not an edge of $G^1$. 
Then $e\cup\{x_1\}$ is a
clique of $G$, a contradiction. Thus $e$ is in 
$E({\rm cl}(G^1))$. Case (ii): Assume $y_1\in e$. Then there is 
$f\in E({\rm cl}(G))$, with $x_1\in f$, such that
$e=(f\setminus\{x_1\})\cup\{y_1\}$. Since $\{x,x_1\}\in E(G)$ for any 
$x$ in $f\setminus\{x_1\}$, one has that $\{x,y_1\}\in E(G^1)$ for 
any $x$ in $f\setminus\{x_1\}$. Then $e$ is a clique of $G^1$. If $e$
is not a maximal clique of $G^1$, there is $x\notin e$ which is
adjacent in $G$ to any vertex of $f\setminus\{x_1\}$ and $x$ is adjacent to
$y_1$ in $G^1$. In particular $x\neq x_1$. 
Then $x$ is adjacent in $G$ to $x_1$ and consequently $x$ is
adjacent in $G$ to any vertex
of $f$, a contradiction because $f$ is a maximal clique of $G$. Thus 
$e$ is in ${\rm cl}(G^1)$. Next we prove the inclusion 
$E({\rm cl}(G^1))\subset E(\mathcal{C}^1)$. Take $e\in E({\rm
cl}(G^1))$, i.e., $e$ is a maximal clique of $G^1$. Case (i): Assume
$y_1\notin e$. Then $e$ is a maximal clique of $G$, and so an edge 
of $\mathcal{C}^1$. Case (ii): Assume $y_1\in e$. Then
$e\setminus\{y_1\}$ is a clique of $G$ and $\{x,y_1\}\in E(G^1)$ for
any $x$ in $e\setminus\{y_1\}$. Then $\{x,x_1\}$ is in $E(G)$ 
for any $x$ in  $e\setminus\{y_1\}$. Hence
$f=(e\setminus\{y_1\})\cup\{x_1\}$ is a clique of $G$. Note that $f$
is a maximal clique of $G$. Indeed if $f$ is not a maximal clique 
of $G$, there is $x\in V(G)\setminus f$ which is adjacent in $G$ to every vertex of
$e\setminus\{y_1\}$ and to $x_1$. Thus $x$ is adjacent to 
$y_1$ in $G^1$ and to every vertex in $e\setminus\{y_1\}$, i.e., 
$e\cup\{x\}$ is a clique of $G^1$, a contradiction. Thus $f\in{\rm
cl}(G)$. Since $e=(f\setminus\{x_1\})\cup\{y_1\}$ we obtain that 
$e\in E(\mathcal{C}^1)$. \QED 

\medskip

Unfortunately we do not have an analogous version 
of Lemma~\ref{dup-comm-cliques} valid for a deletion. In other words,
if $G$ is a graph, the equality ${\rm cl}(G)^w={\rm cl}(G^w)$, with
$w$ an integer vector, fails in general (see
Remark~\ref{del-not-comm-cliques}).

\begin{Remark}\label{del-not-comm-cliques}\rm Let $G$ be a graph. Let $G^1=G\setminus\{x_1\}$  
$($resp. ${\rm cl}(G)^1={\rm cl}(G)\setminus\{x_1\}$$)$ be the graph $($resp.
clutter$)$ obtained from $G$ $($resp. ${\rm cl}(G)$$)$ by deleting the
vertex $x_1$. The equality ${\rm cl}(G)^1={\rm cl}(G^1)$ fails in
general. For instance if $G$ is a cycle of length three, 
then $E({\rm cl}(G)^1)=\emptyset$ and 
${\rm cl}(G^1)$ has exactly one edge. 
\end{Remark}

Let $\mathcal{D}$ be a {\it digraph\/}, that is, $\mathcal{D}$ consists of a
finite set $V(\mathcal{D})$ of vertices and a set $E(\mathcal{D})$ of
ordered pairs of distinct vertices called edges. Let $A$, $B$ be two sets of
vertices of $\mathcal{D}$. For use below recall that a (directed) 
path of $\mathcal{D}$ is called an $A$--$B$ {\it
path} if it runs from a vertex in $A$ to a vertex in $B$. A set $C$ of
vertices is called an $A$--$B$ {\it disconnecting\/} set if $C$
intersects each $A$--$B$ path. For convenience we recall the 
following classical result.

\begin{Theorem}{\rm (Menger's theorem, see
\cite[Theorem~9.1]{Schr2})}\label{menger}\label{menger-theorem} Let $\mathcal{D}$ be a
digraph and let $A$, $B$ be two subsets of $V(\mathcal{D})$. Then the
maximum number of vertex-disjoint $A$--$B$ paths is equal to the
minimum size of an $A$--$B$ disconnecting vertex set.  
\end{Theorem}

We come to the main result of this section.

\begin{Theorem}\label{comp-max-cliques-mfmc} Let $P=(X,\prec)$ be a 
poset on the vertex  
set $X=\{x_1,\ldots,x_n\}$ and let $G$ be its 
comparability graph. If 
$\mathcal{C}={\rm cl}(G)$ is 
the clutter of maximal cliques of $G$, then 
$\mathcal{C}$ satisfies the max-flow min-cut property. 
\end{Theorem}

\demo We can regard $P$ as a transitive digraph without cycles of
length two with vertex set $X$ and edge set $E(P)$, i.e., the edges of
$P$ are ordered pairs $(a,b)$ of 
distinct vertices with $a\prec b$ such that: 

\medskip

(i)\ \,$(a,b)\in E(P)$ and $(b,c)\in E(P)$ $\Rightarrow$ $(a,c)\in
E(P)$ and

(ii) $(a,b)\in E(P)$ $\Rightarrow$ $(b,a)\notin E(P)$.

\medskip

Note that because of these two conditions, 
$P$ is in fact an acyclic digraph, that is, it has no directed 
cycles. 
Let $x_1$ be a vertex of $P$ and let $y_1$ be a new vertex. 
Consider the digraph $P^1$ with vertex set $X^1=X\cup\{y_1\}$ 
and edge set 
$$
E(P^1)=E(P)\cup\{(y_1,x)\vert\, (x_1,x)\in E(P)\}
\cup \{(x,y_1)\vert\, (x,x_1)\in E(P)\}.
$$
The digraph $P^1$ is transitive. Indeed let $(a,b)$ and $(b,c)$ be two
edges of $P^1$. If $y_1\notin\{a,b,c\}$, then $(a,c)\in E(P)\subset
E(P^1)$ because $P$ is transitive. If $y_1=a$, then $(x_1,b)$ and
$(b,c)$ are in $E(P)$. Hence $(x_1,c)\in E(P)$ and $(y_1,c)\in
E(P^1)$. The cases $y_1=b$ and $y_1=c$ are treated similarly. Thus
$P^1$ defines a poset $(X^1,\prec^{1})$. The comparability graph $H$
of $P^1$ is precisely the graph $G^1$ obtained  from $G$ by
duplicating the vertex $x_1$ by the vertex $y_1$. To see 
this note that $\{x,y\}$ is an edge of $G^1$ if and only if $\{x,y\}$
is an edge of $G$ or $y=y_1$ and $\{x,x_1\}$ is an edge of $G$. Thus
$\{x,y\}$ is an edge of $G^1$ if and only if $x$ is related to $y$ in
$P$ or $y=y_1$ and $x$ is related to $y$ in $P^1$, i.e., $\{x,y\}$ is
an edge of $G^1$ if and only if $\{x,y\}$ is an edge of $H$. 
From Lemma~\ref{dup-comm-cliques} we get that 
${\rm cl}(G)^1={\rm cl}(G^1)$, where ${\rm cl}(G)^1$ is the clutter 
obtained from ${\rm cl}(G)$ by duplicating the vertex $x_1$ by the
vertex $y_1$. Altogether we obtain that the clutter ${\rm cl}(G)^1$ 
is the clique clutter of the comparability graph $G^1$ of the poset
$P^1$.

By Theorem~\ref{mfmc-iff-cw} it suffices to prove that 
${\rm cl}(G)^w$ has the K\"onig property 
for all $w\in\mathbb{N}^n$. Since duplications commute with 
deletions, by permuting vertices, we may assume that $w=(w_1,\ldots,w_r,0,\ldots,0)$, 
where $w_i\geq 1$ for $i=1,\ldots,r$. Consider the clutter
$\mathcal{C}_1$ obtained from ${\rm cl}(G)$ by duplicating $w_i-1$ times
the vertex $x_i$ for $i=1,\ldots,r$. We denote the vertex set of
$\mathcal{C}_1$ by $X_1$. By successively applying 
the fact that ${\rm cl}(G)^1={\rm cl}(G^1)$, we conclude that there is a 
poset $P_1$ with comparability graph $G_1$ and vertex set
$X_1$ such that $\mathcal{C}_1={\rm cl}(G_1)$. As before we regard
$P_1$ as a transitive acyclic digraph. 

Let $A$ and $B$ be the set of
minimal and maximal elements of the poset $P_1$, i.e., the elements of $A$ and
$B$ are the sources and sinks of $P_1$ respectively. We set
$S=\{x_{r+1},\ldots,x_n\}$. Consider the digraph
$\mathcal{D}$ whose vertex set is
$V(\mathcal{D})=X_1\setminus S$ and whose edge set is
defined as follows. A pair $(x,y)$ in $V(\mathcal{D})\times
V(\mathcal{D})$ is in $E(\mathcal{D})$ if and
only if $(x,y)\in E(P_1)$ and there is no vertex $z$ in $X_1$ 
with $x\prec z\prec y$. Notice that 
$\mathcal{D}$ is a sub-digraph of $P_1$ 
which is not necessarily the digraph of a
poset. We set $A_1=A\setminus S$ and 
$B_1=B\setminus S$. Note that $\mathcal{C}^w=
\mathcal{C}_1\setminus S$, the clutter obtained from $\mathcal{C}_1$
by removing all vertices of $S$ and all edges sharing a vertex with
$S$. If every edge of $\mathcal{C}_1$ intersects $S$, then
$E(\mathcal{C}^w)=\emptyset$ and there is nothing to prove. Thus we
may assume that there is a maximal clique $K$ of $G_1$ disjoint form
$S$.  Note that by the maximality of $K$ and by the transitivity of
$P_1$ we get that $K$ contains at least one source and one sink of
$P_1$, i.e., $A_1\neq\emptyset$ and $B_1\neq\emptyset$ (see argument
below).

The maximal cliques of $G_1$ not
containing any vertex of $S$ correspond exactly to the 
$A_1$--$B_1$ paths of $\mathcal{D}$. Indeed let
$c=\{v_1,\ldots,v_s\}$ be a maximal clique of $G_1$ disjoint from $S$.
Consider the sub-poset $P_c$ of $P_1$ induced by $c$. Note that $P_c$
is a tournament, i.e., $P_c$ is an oriented graph (no-cycles of length
two) such that any two vertices of $P_c$ are comparable. By
\cite[Theorem~1.4.5]{digraphs} any tournament has a Hamiltonian path,
i.e., a spanning oriented path. Therefore we may assume that
$$ 
v_1\prec v_2\prec\cdots\prec v_{s-1}\prec v_s
$$
By the maximality of $c$ we get that $v_1$ is a source of $P_1$, $v_s$ is a
sink of $P_1$, and $(v_i,v_{i+1})$ is an edge of $\mathcal{D}$ for
$i=1,\ldots,s-1$. Thus $c$ is an $A_1$--$B_1$ path of
$\mathcal{D}$, as required. Conversely let $c=\{v_1,\ldots,v_s\}$ be an
$A_1$--$B_1$ path of $\mathcal{D}$. Clearly $c$ is a clique of $P_1$
because $P_1$ is a poset. Assume that $c$ is not a maximal
clique of $G_1$. Then there is a vertex $v\in X_1\setminus c$ such
that $v$ is related to every vertex of $c$. Since $v_1,v_s$ are a
source and a sink of $P_1$ respectively we get $v_1\prec v\prec v_s$. We claim
that $v_i\prec v$ for $i=1,\ldots,s$. By induction assume that 
$v_i\prec v$ for some $1\leq i<s$. If $v\prec v_{i+1}$, then 
$v_i\prec v\prec v_{i+1}$, a contradiction to the fact that
$(v_i,v_{i+1})$ is an edge of $\mathcal{D}$. Thus $v_{i+1}\prec v$. 
Making $i=s$ we get that $v_s\prec v$, a contradiction. This proves
that $c$ is a maximal clique of $G_1$. Therefore, since the 
maximal cliques of $G_1$ not
containing any vertex in $S$ are exactly the edges of 
$\mathcal{C}^w=\mathcal{C}_1\setminus S$, by 
Menger's theorem (see Theorem~\ref{menger-theorem}) we obtain 
that $\beta_1(\mathcal{C}^w)=\alpha_0(\mathcal{C}^w)$, i.e., 
$\mathcal{C}^w$ satisfies  the K\"onig property. \QED 

\medskip

Let $G$ be a graph. The matrix $A$ whose column vectors are the
characteristic vectors of the maximal cliques of $G$ is called the 
{\it vertex-clique matrix\/} of $G$. It is well known that if $G$ is
a comparability graph 
and $A$ is the
vertex-clique matrix of $G$, then $G$ is
perfect \cite[Corollary~66.2a]{Schr2} and the polytope 
$$
P(A)=\{x\vert\,x\geq 0;\, xA\leq\mathbf{1}\}
$$ 
is integral \cite[Corollary~65.2e]{Schr2}.  The next result
complement this fact. 

\begin{Corollary}\label{compara-int} Let $G$ be a comparability graph
and let $A$ be the vertex-clique matrix of $G$. Then the 
polyhedron $Q(A)=\{x\vert\, x\geq 0;\, xA\geq{\mathbf 1}\}$ 
is integral. 
\end{Corollary}

\demo By Theorem~\ref{comp-max-cliques-mfmc} the clique clutter 
${\rm cl}(G)$ has the max-flow min-cut property. Thus the system 
$xA\geq\mathbf{\mathbf 1}$; $x\geq 0$ is totally dual integral, i.e.,
the maximum in Eq.~(\ref{jun6-2-03-1}) has an integer
optimum solution $y$ for each integer vector $w$ with finite
maximum. Hence $Q(A)$ has only integer vertices by \cite[Theorem
5.22]{Schr2}. \QED

\section{Normally torsion freeness and normality}\label{ntf-normal-sec}

Let $\mathcal{C}$ be a clutter on the vertex set $X$ and 
let $I=I(\mathcal{C})\subset R$ be its edge
ideal. A subset $C\subset X$ is called a 
{\it minimal vertex cover\/} of $\cal C$ if: 
(i) every edge of $\cal C$ contains at least one vertex of $C$, 
and (ii) there is no proper subset of $C$ with the first 
property. Recall that 
$\mathfrak{p}$ is a minimal prime of $I$ if and only if 
$\mathfrak{p}=(C)$ for some minimal vertex cover $C$ of $\mathcal{C}$ 
\cite[Proposition~6.1.16]{monalg}. Thus if 
$C_1,\ldots,C_s$ are the minimal vertex covers of
$\mathcal{C}$, then the primary decomposition of $I$ is
\begin{equation}
I=\mathfrak{p}_1\cap\mathfrak{p}_2\cap\cdots\cap\mathfrak{p}_s,
\end{equation}
where $\mathfrak{p}_i$ is the prime ideal of $R$ generated by $C_i$. The 
$i${\it th} {\it symbolic power} of $I$, denoted by $I^{(i)}$, 
is given by 
$I^{(i)}=\mathfrak{p_1}^i\cap\cdots\cap \mathfrak{p}_s^i$.

\begin{Theorem}[\rm\cite{clutters}]\label{noclu1}
Let $\mathcal{C}$ be a clutter, let $A$ be the incidence matrix of
$\mathcal{C}$, and let $I=I(\mathcal{C})$ be its edge ideal. Then the
following 
are equivalent{\rm :} 
\begin{description}
\item{\rm(i)\ \,} $I$ is normal and $Q(A)=\{x\vert\, x\geq 0;\, xA\geq
\mathbf{1}\}$ is an integral
polyhedron.
\vspace{-1mm}
\item{\rm (ii)\ } $I$ is normally torsion free, 
i.e., $I^{i}=I^{(i)}$ for $i\geq 1$.\vspace{-1mm}
\item{\rm (iii)\,} $\mathcal{C}$ has the max-flow min-cut
property.  
\end{description}
\end{Theorem}

There are some other nice characterizations of the normally
torsion free property that can be found in \cite{reesclu,HuSV}. 

\medskip

Our main algebraic result is:

\begin{Theorem}\label{compara-ntf-normal} If $G$ is a comparability
graph and ${\rm cl}(G)$ is its clique clutter,
then the edge ideal $I=I({\rm cl}(G))$ of ${\rm cl}(G)$ is normally
torsion free and normal.   
\end{Theorem}

\demo It follows from Theorems~\ref{comp-max-cliques-mfmc} and
\ref{noclu1}. \QED

\paragraph{Complete admissible uniform clutters} In this paragraph we
introduce a family of clique clutters of comparability graphs. Let
$d\geq 2$, $g\geq 2$ be two integers and let 
$$
X^1=\{x_1^1,\ldots,x_g^1\},\, X^2=\{x_1^2,\ldots,x_g^2\},\, \ldots, 
X^d=\{x_1^d,\ldots,x_g^d\}
$$
be disjoint sets of variables. The clutter $\mathcal{C}$ with vertex set 
$X=X^1\cup\cdots\cup X^d$ and edge set
$$
E(\mathcal{C})=\{\{x_{i_1}^1,x_{i_2}^2,\ldots,x_{i_d}^d\}\vert\, 
1\leq i_1\leq i_2\leq\cdots \leq i_d\leq g\}
$$
is called a {\it complete admissible uniform clutter}. The edge ideal
of this clutter was introduced and studied in \cite{floystad}. This ideal
has many good properties, for instance $I(\mathcal{C})$ and its Alexander dual
are Cohen-Macaulay and have linear resolutions (see \cite[Proposition
4.5, Lemma 4.6]{floystad}). For a thorough study of Cohen-Macaulay admissible 
clutters see \cite{linearquotients,MRV}. 

\begin{Theorem}\label{cauc-ntf} If $\mathcal{C}$ is a complete admissible uniform
clutter, then its edge ideal $I(\mathcal{C})$ is normally torsion
free and normal. 
\end{Theorem}

\demo Let $P=(X,\prec)$ be the poset with vertex set $X$ and partial
order given by $x_k^\ell\prec x_p^m$ if and only if $1\leq \ell<m\leq
d$ and $1\leq k\leq p\leq g$. We denote the comparability graph of $P$
by $G$. We claim that $E(\mathcal{C})=E({\rm cl}(G))$, where ${\rm
cl}(G)$ is the clique clutter of $G$. Let
$f=\{x_{i_1}^1,x_{i_2}^2,\ldots,x_{i_d}^d\}$ be an edge of
$\mathcal{C}$, i.e.,
we have $1\leq i_1\leq i_2\leq\cdots \leq i_d\leq g$. Clearly $f$ is a
clique of $G$. If $f$ is not maximal, then there is a vertex
$x_k^\ell$ not in $f$ which is adjacent in $G$ to every vertex 
of $f$. In particular $x_k^\ell$ must be comparable to
$x_{i_\ell}^\ell$, which is impossible. Thus $f$ is an edge 
of  ${\rm cl}(G)$. Conversely let $f$ be an edge of ${\rm cl}(G)$. 
We can write $f=\{x_{i_1}^{k_1},x_{i_2}^{k_2},\ldots,x_{i_s}^{k_s}\}$,
where $k_1<\cdots< k_s$ and $i_1\leq\cdots\leq i_s$. By the maximality
of $f$ we get that $s=d$ and $k_i=i$ for $i=1,\ldots,d$. Thus 
$f$ is an edge of $\mathcal{C}$. Hence by 
Theorem~\ref{compara-ntf-normal} we
obtain that $I(\mathcal{C})$ is normally torsion free and normal. \QED

\medskip
\noindent
{\bf Acknowledgments.} We thank Seth Sullivant
for pointing out an alternative  
proof of Theorem~\ref{compara-ntf-normal} based on the fact that the
edge ideal of a comparability 
graph is differentially perfect (see
\cite[Section~4]{sullivant-csp}).  One of the 
consequences of being differentially perfect is that if every maximal
chain in the poset has the same fixed length $k$, then the edge ideal
of the clique clutter of its comparability graph is
normally torsion free. The general case, that is, 
the situation where not all maximal chains have the same length, can
be reduced to the special case above.  We also thank Alexander Schrijver for pointing out 
\cite[Theorem~14.18]{Schr2} and the characterization of the max-flow
min-cut property given in Theorem~\ref{mfmc-iff-cw}. Finally, we thank the referee for a
careful reading of the paper and for pointing out \cite{floystad} as
an earlier reference for complete admissible uniform clutters. 
The authors have been informed 
that Corollary~\ref{crit-rounding-normali} was observed by Ng\^{o}
Vi\^{e}t Trung 
when $I$ is the edge ideal of a hypergraph.

\bibliographystyle{plain}

\begin{thebibliography}{10}

\bibitem{digraphs} J. Bang-Jensen and G. Gutin, {\it Digraphs}. {\it Theory,
Algorithms and Applications\/}, Springer Monographs in Mathematics,
Springer, 2nd printing, 2006.

\bibitem{baum-trotter} S. Baum and L. E. Trotter, 
Integer rounding for polymatroid and branching optimization problems,
SIAM J. Algebraic Discrete Methods {\bf 2} (1981), no. 4, 416--425. 

\bibitem{BHer}{W. Bruns and J. Herzog, 
{\em Cohen-Macaulay Rings\/},  Cambridge University 
Press, Cambridge, Revised Edition, 1997.}


\bibitem{normali} C. Escobar, R. H. Villarreal and Y. Yoshino, Torsion
freeness and normality of blowup rings of monomial ideals, 
{\it Commutative Algebra\/}, Lect. Notes Pure Appl. Math. 
{\bf 244}, Chapman \& Hall/CRC, Boca Raton, FL, 2006, pp. 69-84. 

\bibitem{floystad} G. F{l\o}ystad and J. E. Vatne, 
(Bi-)Cohen-Macaulay simplicial complexes and their associated
coherent sheaves, Comm. Algebra {\bf 33} (2005), no. 9, 
3121--3136. 

\bibitem{reesclu}{I. Gitler, E. Reyes and R. H. Villarreal, 
Blowup algebras of square--free monomial ideals and some links to
combinatorial optimization problems, Rocky Mountain J. Math., to appear.} 


\bibitem{clutters}{I. Gitler, C. Valencia and R. H. Villarreal, 
A note on Rees algebras and the MFMC property, Beitr\"age Algebra Geom. {\bf 48}
(2007), No. 1, 141-150.}

\bibitem{linearquotients} H. T. H\`{a}, S. Morey and R.H. Villarreal,
Cohen-Macaulay 
admissible clutters, Journal of Commutative Algebra, to appear.

\bibitem{HuSV}{C. Huneke, A. Simis and W. V. Vasconcelos, Reduced
normal cones are domains, Contemp.  Math. {\bf 88} (1989),
95--101.}

\bibitem{MRV} S. Morey, E. Reyes and R. H. Villarreal, 
Cohen-Macaulay, shellable and unmixed clutters with a perfect
matching of K\"{o}nig type, J. Pure Appl. 
Algebra {\bf 212} (2008), 1770-1786.

\bibitem{Schr}{A. Schrijver, {\it Theory of Linear and 
Integer Programming\/}, John Wiley \& Sons, New York, 1986.}

\bibitem{Schr2} {A. Schrijver,{\it Combinatorial Optimization\/}, 
Algorithms and Combinatorics {\bf 24}, Springer-Verlag, Berlin, 2003.}

\bibitem{sullivant-csp} 
S. Sullivant, Combinatorial symbolic powers, J. Algebra
{\bf 319(1)} (2008), 115-142. 

\bibitem{Vas1}{W. V. Vasconcelos, {\it Computational Methods in
Commutative Algebra and Algebraic Geometry\/}, 
Springer-Verlag, 1998.}

\bibitem{bookthree} W. V. Vasconcelos, {\it Integral Closure\/}, Springer
Monographs in Mathematics, Springer, New York, 2005.  

\bibitem{monalg}{R. H. Villarreal, {\it Monomial Algebras\/}, 
Monographs and 
Textbooks in Pure and Applied Mathematics {\bf 238}, Marcel 
Dekker, New York, 2001.} 

\end{thebibliography}

\end{document}